\input amstex\documentstyle{amsppt}  
\pagewidth{12.5cm}\pageheight{19cm}\magnification\magstep1
\topmatter
\title Strata and almost special representations\endtitle
\author G. Lusztig\endauthor
\address{Department of Mathematics, M.I.T., Cambridge, MA 02139}\endaddress
\thanks{Supported by NSF grant DMS-2153741}\endthanks
\endtopmatter   
\document
\define\mat{\matrix}
\define\endmat{\endmatrix}

\define\sg{\text{sgn}}

\define\Irr{\text{\rm Irr}}

\define\mpb{\medpagebreak}

\define\si{\sim}

\define\sqc{\sqcup}

\define\part{\partial}

\define\m{\mapsto}
\define\do{\dots}

\define\lra{\leftrightarrow}

\define\sub{\subset}

\define\nl{\newline}
\redefine\i{^{-1}}

\define\ot{\otimes}

\define\ph{\phi}

\define\s{\sigma}

\define\th{\theta}

\redefine\l{\lambda}

\redefine\G{\Gamma}
\redefine\D{\Delta}

\define\Si{\Sigma}
\define\Th{\Theta}
\redefine\L{\Lambda}
\define\Ph{\Phi}

\define\kk{\bold k}

\define\CC{\bold C}

\define\QQ{\bold Q}

\define\ca{\Cal A}

\define\cs{\Cal S}

\define\cu{\Cal U}

\define\sha{\sharp}

\head Introduction\endhead
\subhead 0.1\endsubhead
Let $G$ be a connected reductive group over an algebraically
closed field $\kk$; let $W$ be the Weyl group of $G$.
Let $\Irr(W)$ be the set of irreducible representations (up to
isomorphism) of $W$ (over $\QQ$).

In \cite{L15} we have defined a partition
$G=\sqc_{\Si\in Str_G}\Si$ of $G$
into finitely many subsets $\Si$ called strata. (Each stratum is a union
of conjugacy classes in $G$ of fixed dimension.) Now $Str_G$ can be
regarded as an enlargement of the set $\cu_G$ of unipotent classes of
$G$ and we expect that various objects attached to
$\cu_G$ have analogues for $Str_G$. One such object is the collection
of the groups of components of the centralizers
(in the adjoint group of $G$) of elements in the various unipotent
classes. In this paper we shall attach to
any stratum $\Si$ of $G$ a finite group $\G_\Si$
(defined up to isomorphism); this attachment seems to be
similar to the attachment in the previous sentence (see 3.3-3.5).
To attach the finite group $\G_\Si$ to a stratum $\Si$
we first attach to $\Si$
(following \cite{L15}) an irreducible representation of $W$,
and then we use the indexing of $\Irr(W)$ by pairs of finite
groups given in \cite{L19}. If $W$ is irreducible, $\G_\Si$
is a product of copies of $S_2$ or a single $S_3,S_4$ or $S_5$
(here $S_n$ is the symmetric group in $n$ letters.)
It would be desirable to have a more
direct definition of this attachment.

\subhead 0.2\endsubhead
In \S1 we show that if $W$ is of simply laced type, 
the representation of $W$ attached to a stratum in \cite{L15} is
almost special in the sense of \cite{L24a}. In \S2 we study the
action of tensoring by sign on
almost special representations. In \S3 we define the finite
group $\G_\Si$ attached to a stratum $\Si$.

\head 1. From $2$-special representations to almost special
representations\endhead
\subhead 1.1\endsubhead
In \cite{L15} the set
$Str_G$ of strata of $G$ was put in bijection $\Si\lra E_\Si$
with a subset $\cs_2(W)$ of $\Irr(W)$; the set
$\cs_2(W)$ depends only on $W$, not on $G$ (its
elements are said to be {\it $2$-special representations} of $W$).

Let $\ph(W)$ be the set of families of $W$; we have
$\Irr(W)=\sqc_{c\in\ph(W)}c$. In \cite{L84} we have attached to any
$c\in\ph(W)$ a certain finite group $\G_c$ and in \cite{L24a} we have
defined a subset $\ca_{\G_c}$ of $c$. Let
$\ca(W)=\sqc_c\ca_{\G_c}\sub\Irr(W)$;
the elements of $\ca(W)$ are said to be the {\it almost special
representations} of $W$.

\proclaim{Theorem 1.2} If $W$ is of simply laced type, we have
$\cs_2(W)\sub\ca(W)$.
\endproclaim
Without the simplylacedness assumption, this does not hold. Indeed, if
$W$ is of type $B_2$, we have $\cs_2(W)=\Irr(W)$, $\ca(W)=\Irr(W)-\{\s\}$
for some one dimensional $\s\in\Irr(W)$.

\subhead 1.3\endsubhead
We can assume that $W$ is irreducible.
If $W$ is of type $A$ we have $\ca_2(W)=\ca(W)=\Irr(W)$ hence the
theorem holds.

Assume now that $W$ is of type $D_n,n\ge4$.
Let $Sy_{n,m}$ be the set of all symbols
$$\L=\left(\mat \l_1 &\l_2&\do&\l_m\\ \mu_1 &\mu_2&\do&\mu_m\endmat
\right)$$
where $\l_i,\mu_i$ are integers such that
$0\le\l_1<\l_2<\do<\l_m,0\le\mu_1<\mu_2<\do<\mu_m$,
$\sum_i\l_i+\sum_i\mu_i=n+m^2-m$.
Consider the equivalence relation on $\sqc_m Sy_{n,m}$ generated by the relations
$$\left(\mat \l_1 &\l_2&\do&\l_m\\ \mu_1 &\mu_2&\do&\mu_m\endmat
\right)\si
\left(\mat 0&\l_1+1 &\l_2+1&\do&\l_m+1\\ 0&\mu_1+1 &\mu_2+1&\do&
\mu_m+1\endmat\right)$$
and 
$$\left(\mat \l_1 &\l_2&\do&\l_m\\ \mu_1 &\mu_2&\do&\mu_m\endmat
\right)\si
\left(\mat \mu_1 &\mu_2&\do&\mu_m\\ \l_1 &\l_2&\do&\l_m\endmat
\right).$$
Let $Sy_n$ be the set of equivalence classes on 
$\sqc_m Sy_{n,m}$. We say that
$$\left(\mat \l_1 &\l_2&\do&\l_m\\ \mu_1 &\mu_2&\do&\mu_m\endmat
\right)\in Sy_{n,m}$$
is degenerate if $\l_1=\mu_1,\l_2=\mu_2,\do,\l_m=\mu_m$.
Let $Sy''_n$ be the set of all elements of $Sy_n$ that are represented
by degenerate symbols; let $Sy'_n=Sy_n-Sy''_n$.

Let $Y_n$ be the subset of $Sy'_n$ consisting of all equivalence
classes of symbols
$$\left(\mat \l_1 &\l_2&\do&\l_m\\ \mu_1 &\mu_2&\do&\mu_m\endmat
\right)\in Sy_{n,m}$$
such that $\l_1\le\mu_1,\l_2\le\mu_2,\do,\l_m\le\mu_m$ (with at least
one of the last inequalities being strict).

Let $X_n$ be the subset of $Sy'_n$ consisting of all equivalence
classes of symbols
$$\left(\mat \l_1 &\l_2&\do&\l_m\\ \mu_1 &\mu_2&\do&\mu_m\endmat
\right)\in Sy_{n,m}$$
such that $\l_1\le\mu_1,\l_2\le\mu_2,\do,\l_m\le\mu_m$ (with at least
one of the last inequalities being strict) and such that
$\mu_1\le\l_2+3,\mu_2\le\l_3+3,\do,\mu_{m-1}\le\l_m+3$.

As in \cite{L84, 4.6} we identify $\Irr(W)$ with the disjoint union of
$Sy'_n$ with two copies of $Sy''_n$.
From \cite{L15, 3.10(b)} (resp. from \cite{L24a}) one can deduce that
under this identification the subset $\cs_2(W)$ (resp. $\ca(W)$) of
$\Irr(W)$ corresponds to the disjoint union of the subset $X_n$ (resp.
$Y_n$) of $Sy'_n$ with two copies of $Sy''_n$.
Since $X_n\sub Y_n$, we see that the theorem holds in our case.

\subhead 1.4\endsubhead
We now assume that $W$ is of type $E_6,E_7$ or $E_8$.
In the following three subsections we give
in each of these cases a table of the objects in $\ca(W)$
(following \cite{L24a}); the objects of $\cs_2(W)$ (known from
\cite{L15}) are exactly the objects in $\ca(W)$ that are not marked
as $((?))$. (The notation for the objects of $\Irr(W)$ is the same as
that of \cite{L15, 4.4}.) In particular we have $\cs_2(W)\sub\ca(W)$. In fact,
$\ca(W)-\cs_2(W)$ consists of $0,1$ or $2$ objects when $W$ is of type
$E_6,E_7$ or $E_8$ respectively. This completes the proof of Theorem 1.2.

In each table of 1.5-1.7 two representations are written in the same row
if and only if they are in the same family of $W$.

\subhead 1.5. Type $E_6$\endsubhead

$1_0$

$6_1$

$20_2$

$30_3,15_4$

$64_4$

$60_5$

$24_6$

$81_6$

$80_7,60_8,10_9$

$81_{10}$

$60_{11}$

$24_{12}$

$64_{13}$

$30_{15},15_{16}$

$20_{20}$

$6_{25}$

$1_{36}$

\subhead 1.6. Type $E_7$\endsubhead

$1_0$

$7_1$

$27_2$

$21_3$

$56_3,35_4$

$120_4,((15_7))$

$189_5$

$315_7,280_8,70_9$

$405_8,216_9$

$378_9$

$210_{10}$

$420_{10},84_{12}$

$105_{12}$

$210_6$

$105_6$

$168_6$

$189_7$

$512_{11}$

$210_{13}$

$420_{13},84_{15}$

$378_{14}$

$105_{15}$

$405_{15},216_{16}$

$315_{16},280_{17},70_{18}$

$189_{20}$

$168_{21}$

$105_{21}$

$210_{21}$

$189_{22}$

$120_{24},15_{28}$ 

$56_{30},35_{31}$

$21_{36}$

$27_{37}$

$7_{46}$

$1_{63}$

\subhead 1.7. Type $E_8$\endsubhead

$1_0$

$8_1$

$35_2$

$112_3,84_4$

$210_4,((50_8))$

$560_5$

$567_6$

$700_6,400_7$

$1400_7,1344_8,448_9$

$1400_8,1050_{10},175_{12}$

$3240_9$

$2268_{10},972_{12}$

$2240_{10},1400_{11}$

$4096_{11}$

$525_{12}$

$4200_{12},840_{14}$

$2800_{13},((700_{16}))$

$4536_{13}$

$2835_{14}$

$6075_{14}$

$4200_{15}$

$5600_{15},3200_{16}$

$4480_{16},7168_{17},4200_{18},3150_{18},2016_{19},
1344_{19},420_{20},168_{24}$

$2100_{20}$

$4200_{21}$

$5600_{21},3200_{22}$

$6075_{22}$

$2835_{22}$

$4536_{23}$

$4200_{24},840_{26}$

$2800_{25},700_{28}$

$4096_{26}$

$2240_{28},1400_{29}$

$2268_{30},972_{32}$

$3240_{31}$

$1400_{32},1050_{34},175_{36}$

$525_{36}$

$1400_{37},1344_{38},448_{39}$

$700_{42},400_{43}$

$567_{46}$

$560_{47}$

$210_{52},50_{56}$    

$112_{63},84_{64}$

$35_{74}$

$8_{91}$

$1_{120}$

\head 2. Tensoring by the sign representation\endhead
\subhead 2.1\endsubhead
In this section we assume that $W$ (as in 1.2) is irreducible. We show:

(a) {\it If $E\in\ca(W)$ then $E\ot\sg\in\ca(W)$ except if $W$ is of
type $E_7$ and $\dim E=512$ or if $W$ is of type $E_8$ and
$\dim E=4096$.}
\nl
If $W$ is of type $A$ there is nothing to prove.
If $W$ is of type $E_6,E_7$ or $E_8$ the result follows from the tables
in 1.5-1.7.

In the remainder of this subsection we assume that $W$ is
of type $D_n,n\ge4$. Let
$$\L=\left(\mat \l_1 &\l_2&\do&\l_m\\ \mu_1 &\mu_2&\do&\mu_m\endmat
\right)\in Sy_{n,m}$$
be such that $\l_1\le\mu_1,\l_2\le\mu_2,\do,\l_m\le\mu_m$.
Let $N$ be an integer such that $N\ge\l_i,N\ge\mu_i$ for all $i$.
We arrange the numbers in $\{0,1,\do,N\}-\{\l_1,\l_2,\do,\l_m\}$
in increasing order $A_1<A_2<\do<A_{N-m+1}$.
We arrange the numbers in $\{0,1,\do,N\}-\{\mu_1,\mu_2,\do,\mu_m\}$
in increasing order $B_1<B_2<\do<B_{N-m+1}$.

We show that 

(b) $B_h\le A_h$
\nl
for any $h=1,2,\do,N-m+1$.

Let $e$ be the largest number in $\{1,2,...,m\}$ such that
$\l_1<A_h,\l_2<A_h,\do,\l_e<A_h$.
We have
$$\align&h=\l_1+(\l_2-\l_1-1)+...+(\l_e-\l_{e-1}-1)
+(A_h-\l_e-1)+1\\&=A_h-\sha(k\in[1,m];\l_k<A_h)-1.\endalign$$
Similarly we have
$$h=B_h-\sha(k\in[1,m];\l_k<B_h)-1.$$
It follows that
$$A_h-\sha(k\in[1,m];\l_k<A_h)=B_h-\sha(k\in[1,m];\l_k<B_h).$$
Assume that $B_h=A_h+x$ for some $x>0$. 
Then for some $z$ we have
$$\sha(k;\l_k<A_h)=z, \sha(k;\mu_k<B_h)=z+x.$$
We have
$$\l_1<A_h,\l_2<A_h,\do,\l_z<A_h,\l_{z+1}>A_h,\do,\l_{z+x}>A_h,$$
$$\mu_1<B_h,\mu_2<B_h,\do,\mu_z<B_h,\mu_{z+1}<B_h,\do,\mu_{z+x}<B_h.$$
Thus for $m\in\{1,2,\do,x\}$ we have $A_h<\l_{z+m}\le\mu_{z+m}<B_h$.
We see that there are at least $x$ distinct integers $y$ such that
$A_h<y<B_h$. This contradicts $B_h=A_h+x$. We see that (b) holds.

For $i=1,2,\do,N-m+1$ we set $A'_i=N-A_{N-m+2-i}$, $B'_i=N-B_{N-m+2-i}$.
From (b) we have $A'_h\le B'_h$ for any $h$.
But it is known (see \cite{L84}) that
$$\L\m\left(\mat A'_1 &A'_2&\do&A'_{N-m+1}\\ B'_1 &B'_2&\do&B'_{N-m+1}
\endmat\right)$$
corresponds to tensoring by $\sg$ in $\Irr(W)$. This proves (a).

\mpb

For example,
$$\left(\mat 0&1&2&3&4\\2&3&4&5&6\endmat\right),$$
$$\left(\mat 0&1\\ 5&6\endmat\right)$$
represent objects of $\ca(W), n=10,$ related to each other by $\ot\sg$.
The first one is in $\cs_2(W)$, the second one is not.

\subhead 2.2\endsubhead
In this subsection we assume that $W$ is of type $E_7$ or $E_8$.
Let $E\in\ca(W)$ be such that $E\ot\sg\in\ca(W)$ (see 2.1).
From the table in 1.6 we see that
If $E\in\cs_2(W)$ then $E\ot\sg\in\cs_2(W)$
except in the following three cases:

(i) $W$ is of type $E_7$ and $E=15_{28}$;

(ii) $W$ is of type $E_8$ and $E=50_{56}$;

(iii) $W$ is of type $E_8$ and $E=700_{28}$.

Let $cl(W)$ be the set of conjugacy classes in $W$ and let
${}'\Ph:cl(W)@>>>\cs_2(W)$ be the surjective map defined in
\cite{L15}. One can show that $E$ in (i),(ii) are characterized
by the property that ${}'\Ph\i(E)$ has a maximum number of elements
(that is $4$ in type $E_7$ and $5$ in type $E_8$) or alternatively by
the property that ${}'\Ph\i(E)$ contains the conjugacy class consisting
of $-1$. Now $E$ in (iii) can be characterized by the property that it is
obtained by truncated induction from the $E$ in (i).

\head 3. Strata and the group $H'/H$\endhead
\subhead 3.1\endsubhead
For any $c\in\ph(W)$ we set $G(c)=\cup_{\Si\in Str_G;E_\Si\in c}\Si$.
We have a partition $G=\sqc_{c\in\ph(W)}G(c)$.

\subhead 3.2\endsubhead
We fix $c\in\ph(W)$.
In \cite{L19, 0.4} (see also \cite{L24})
we have defined a bijection $\th:\Th_c@>\si>>c$
where $\Th_c$ is a certain set of pairs $(H,H')$ of subgroups of $\G_c$
defined up to simultaneous
$\G_c$-conjugacy with $H$ being a normal subgroup of $H'$.
Let $\Th_c^{\cs_2}=\th\i(\cs_2\cap c)$; then $\th$ restricts to a bijection
$\Th_c^{\cs_2}@>\si>>\cs_2\cap c$.
We now attach to any stratum $\Si$ of $G$ contained in
$G(c)$ the finite group $\G_\Si=H'/H$ where $(H,H')\in\Th_c^{\cs_2}$
corresponds to the representation in $\cs_2(W)\cap c$ attached to the
stratum.

For example, if $\Si$ is the stratum in $G_(c)$ that contains
a special unipotent class then $\G_\Si=\G_c/\{1\}=\G_c$.

\subhead 3.3\endsubhead
In this subsection we assume that
$G$ is of type $E_8$ and $c\in\ph(W)$ is such that $\G_c=S_5$.

We list the strata $\Si$ in $G(c)$ (in terms
of $E_\Si\in c$) and the corresponding quotient $\G_\Si=H'/H$.

$4480_{16}....S_5/1=S_5$

$7168_{17}.... S_3S_2/S_2=S_3$

$4200_{18}.... \D_8/S_2S_2=S_2$ 

$3150_{18}....S_3S_2/S_3=S_2$

$2016_{19}.... S_3S_2/S_3S_2=1$

$1344_{19}.... S_4/S_4=1$

$420_{20}....S_5/S_5=1$

$168_{24}.... \D_8/\D_8=1$.
\nl
(We use the obvious notation for the subgroups of $S_5$; $\D_8$ is the
dihedral group of order $8$.)

We now assume that $\kk=\CC$. In the first seven cases, $\Si$
contains a unipotent class and the group $\G_\Si$ is the group of
components of the centralizer of an element in that class. In the eighth
case $\Si$ does not contain a unipotent class.

\subhead 3.4\endsubhead
In this subsection we assume that
$G$ is of type $F_4$ and $c\in\ph(W)$ is such that $\G_c=S_4$.

We list the strata $\Si$ in $G(c)$ (in terms
of $E_\Si\in c$) and the corresponding quotient $\G_\Si=H'/H$.

$12_4....S_4/1=S_4$

$16_5.... S_2S_2/S_2=S_2$

$9_6....\D_8/S_2S_2=S_2$

$6_6.... S_3/S_3=1$ 

$4_7....S_2S_2/S_2=S_2$

$9_6....S_3/1=S_3$

$4_7.... S_4/S_4=1$

$4_8.... S_2S_2/S_2S_2=1$
\nl
(We use the obvious notation for the subgroups of $S_4$; $\D_8$ is the
dihedral group of order $8$.)

We now assume that $\kk=\CC$. In the first five cases, $\Si$
contains a unipotent class and the group $\G_\Si$ is the group of
components of the centralizer of an element in that class. In the remaining cases $\Si$ does not contain a unipotent class.

\subhead 3.5\endsubhead
In this subsection we assume that
$G$ is of type $G_2$ and $c\in\ph(W)$ is such that $\G_c=S_3$.

We list the strata $\Si$ in $G_c$ (in terms
of $E_\Si\in c$) and the corresponding quotient $\G_\Si=H'/H$.

$2_1....S_3/1=S_3$

$2_2....S_2/S_2=1$

$1_3....S_3/S_3=1$

$1_3....S_2/1=S_2$
\nl
(We use the obvious notation for the subgroups of $S_3$.)

We now assume that $\kk=\CC$. In the first three cases, $\Si$
contains a unipotent class and the group $\G_\Si$ is the group of
components of the centralizer of an element in that class. In the
remaining case $\Si$ does not contain a unipotent class.

\enddocument